\documentclass[a4paper,11pt]{article}
\usepackage{latexsym}
\usepackage{amsmath}
\usepackage{amssymb}
\usepackage{enumerate}
\usepackage{theorem}
\usepackage{array}
\newtheorem{theorem}{Theorem}[section]
\newtheorem{lemma}[theorem]{Lemma}
\newtheorem{corollary}[theorem]{Corollary}

\theorembodyfont{\rmfamily}

\newtheorem{setup}[theorem]{Setup}
\newtheorem{remark}[theorem]{Remark}

\newcommand{\proof}{\noindent \mbox{\em Proof.\hspace*{2mm}}}
\newcommand{\qed}{\hfill \mbox{$  \Box $}}

\title{Filling curves for $\mathbb{P}^1 \times \mathbb{P}^1$}
\author{
Masaaki Homma
\\
 Department of Mathematics and Physics\\
Kanagawa University\\
Hiratsuka 259-1293, Japan\\
homma@kanagawa-u.ac.jp
\and
Seon Jeong Kim
\thanks{Partially 
supported by Basic Science Research Program through the National Research Foundation of
Korea (NRF) funded by the Ministry of Education (2016R1D1A1B01011730).}
\\
 Department of Mathematics, and RINS\\
Gyeongsang National University\\
Jinju 52828, Korea \\
skim@gnu.kr
}
\date{}

\begin{document}
\maketitle
\begin{abstract}
We determine the minimal bi-degree(s) of an irreducible filling curve over $\mathbb{F}_q$ for $\mathbb{P}^1\times \mathbb{P}^1$.
It is $(q+1, q+1)$ if $q\neq 2$, and they are $(4,3)$ and $(3,4)$ if $q=2$.
\end{abstract}

\section{Introduction}
A plane filling curve over $\mathbb{F}_q$ is a curve $C$ in $\mathbb{P}^2$ such that the set of $\mathbb{F}_q$-points $C(\mathbb{F}_q)$ of $C$ 
coincides with $\mathbb{P}^2(\mathbb{F}_q)$ of $\mathbb{P}^2$.
If $C$ is a plane filling curve of degree $d$,
then $d\geq q+1$ and equality occurs only if $C$ is a union of $q+1$ $\mathbb{F}_q$-lines passing through a common $\mathbb{F}_q$-point.
Hence if a plane filling curve is irreducible, then $d \geq q+2$;
and such a curve of degree $q+2$ actually exists, which is due to
Tallini \cite{tal2}.
The family of plane filling curves of degree $q+2$ has been studied
 by several authors \cite{dur2018, hom2020, hom-kim2013, tal1, tal2}
from several points of view, for example, the specific defining equations,
classification, their automorphism group, and their ``degenerations".

We want to study curves on $\mathbb{P}^1 \times \mathbb{P}^1$ from a similar point of view.
As the first step, in this paper, we determine the minimal bi-degree of irreducible filling curves for $\mathbb{P}^1 \times \mathbb{P}^1$.

Let $C$ be an irreducible curve on  $\mathbb{P}^1 \times \mathbb{P}^1$ defined by
a bi-homogeneous polynomial
$
F(X_0, X_1; Y_0, Y_1) \in \mathbb{F}_q[X_0, X_1, Y_0, Y_1]
$
of bi-degree $(a,b)$.
Here the first two variables $X_0, X_1$ are coordinates of the first factor $\mathbb{P}^1$ of  $\mathbb{P}^1 \times \mathbb{P}^1$ and the second $Y_0, Y_1$ are those of the second one.

It is easy to show that if
\begin{equation}\label{filling}
C(\mathbb{F}_q) = \mathbb{P}^1(\mathbb{F}_q) \times \mathbb{P}^1(\mathbb{F}_q),
\end{equation}
then $a \geq q+1$ and $b \geq q+1$
(Lemma~\ref{mainlemma}).

If a curve $C \subset \mathbb{P}^1 \times \mathbb{P}^1$ has the property (\ref{filling}), the curve is called a $\mathbb{P}^1 \times \mathbb{P}^1$-filling curve over $\mathbb{F}_q$.

For the existence problem,
we can prove that
\begin{enumerate}[(i)]
\item when $q \neq 2$, there is an irreducible  $\mathbb{P}^1 \times \mathbb{P}^1$-filling curve  over $\mathbb{F}_q$ of bi-degree $(q+1, q+1)$;
\item when $q=2$, there are no irreducible  $\mathbb{P}^1 \times \mathbb{P}^1$-filling curves  over $\mathbb{F}_q$ of bi-degree $(3, 3)$,
but there are irreducible  $\mathbb{P}^1 \times \mathbb{P}^1$-filling curves  over $\mathbb{F}_q$ of bi-degree $(4, 3)$ and $(3,4)$.
\end{enumerate}

At the end of Introduction, we add a parenthetic remark.
For an algebraic variety $X$ over $\mathbb{F}_q$,
$N_q(X)$ denotes the cardinality of the set of $\mathbb{F}_q$-points of $X$.
Let $C$ be a nondegenerate irreducible curve of degree $d$ over $\mathbb{F}_q$ in $\mathbb{P}^r$.
Then
\[
N_q(C) \leq
\left\lfloor\frac{(q-1)(q^{r+1}-1)}{q(q^r-1) -r(q-1)}d \right\rfloor
\]
which was shown in \cite[Theorem~3.2]{hom2012}.
The curve of bi-degree $(4,3)$ (and also $(3,4)$) in (ii) above
attains this upper bound
via Segre embedding $\mathbb{P}^1 \times \mathbb{P}^1 \hookrightarrow \mathbb{P}^3$.
Actually,
\[
\left\lfloor\frac{(q-1)(q^{r+1}-1)}{q(q^r-1) -r(q-1)}d \right\rfloor
= \left\lfloor\frac{15}{2\cdot 7 -3}7 \right\rfloor =\lfloor 9.5...\rfloor =9
\]
which coincides with $N_2(\mathbb{P}^1 \times \mathbb{P}^1)$.
\section{A lemma}

\begin{lemma}\label{mainlemma}
Let $F(X_0,X_1; Y_0, Y_1) \in \mathbb{F}_q[X_0,X_1; Y_0, Y_1]$
be an irreducible bi-homogeneous polynomial of bi-degree $(a,b)$,
and $C$ be the curve on $\mathbb{P}^1 \times \mathbb{P}^1$
defined by $F=0$.
Suppose that
\[
C(\mathbb{F}_q) = \mathbb{P}^1(\mathbb{F}_q) \times \mathbb{P}^1(\mathbb{F}_q).
\]
Then $a \geq q+1$ and $b \geq q+1$.
More precisely, there are bi-homogeneous polynomials $f$ and $g$
in $\mathbb{F}_q[X_0,X_1; Y_0, Y_1]$ of bi-degree $(a -(q+1), b)$
and $(a, b -(q+1)) $ respectively such that
\[
F = f \cdot (X_0^q X_1 - X_0 X_1^q) + g \cdot(Y_0^q Y_1 - Y_0 Y_1^q).
\]
\end{lemma}
\proof
Let $\pi_1: \mathbb{P}^1 \times \mathbb{P}^1 \to \mathbb{P}^1$
be the first projection.
Since $C$ is irreducible, $\pi_1^{-1}(P)$ is not a component of $C$
for any $P \in \mathbb{P}^1(\mathbb{F}_q)$.
Hence
\[
N_q(C) \leq \sum_{P \in \mathbb{P}^1(\mathbb{F}_q)} (C.\pi_1^{-1}(P)) =b(q+1).
\]
Since $N_q(C)= (q+1)^2$,
we have $b \geq q+1$.
Similarly, $ a \geq q+1$.

Next we show the additional assertion.
Let $x = X_1/X_0$ and $y= Y_1/Y_0$.
Then
\[
F(1,x;1,y) = \frac{F(X_0,X_1; Y_0, Y_1)}{X_0^a Y_0^b}
\]
is a polynomial in $x$ and $y$,
and vanishes on
\[
\mathbb{A}^2(\mathbb{F}_q) =
 \left(
(\mathbb{P}^1 \setminus \{X_0 =0\}) \times 
(\mathbb{P}^1  \setminus \{Y_0 =0\})
 \right)(\mathbb{F}_q).
\]
Hence, there are polynomials
$u(x,y), v(x,y) \in \mathbb{F}_q[x,y]$
such that
\[
F(1,x;1,y) =u(x,y)(x -x^q) + v(x,y)(y-y^q).
\]
Therefore, there are bi-homogeneous polynomials
$U=U(X_0, X_1;Y_0,Y_1)$ and $V=V(X_0, X_1;Y_0,Y_1)$
of bi-degree $(a-q, b)$ and $(a, b-q)$ respectively
such that
\begin{equation}\label{eq1}
F = U\cdot (X_0^{q-1}X_1 - X_1^q) + V\cdot (Y_0^{q-1}Y_1-Y_1^q).
\end{equation}
Then
\begin{equation}\label{eq2}
F(0, X_1; Y_0,Y_1)
= -U(0,X_1;Y_0,Y_1) X_1^q +V(0,X_1; Y_0, Y_1) (Y_0^{q-1}Y_1-Y_1^q).
\end{equation}
Since $U$ is of bi-degree $(a-q, b)$,
\begin{equation}\label{eq3}
U(0,X_1;Y_0,Y_1)= f_1(Y_0, Y_1)X_1^{a-q}
\end{equation}
for a certain homogeneous polynomial $f_1(Y_0, Y_1)$ in $Y_0, Y_1$
of degree $b$.
Since the leftmost  and the rightmost terms of (\ref{eq2})
vanish on $\{ (0,1) \times (1, \beta ) \mid \beta \in \mathbb{F}_q\}$,
$f_1(1, \beta)=0$ 
for any $\beta \in \mathbb{F}_q$.
Hence there is a homogeneous polynomial $f_2(Y_0, Y_1)$ of degree $b-q$
such that
\begin{equation}\label{eq4}
f_1(Y_0, Y_1) = f_2(Y_0, Y_1)(Y_0^{q-1} Y_1 - Y_1^q).
\end{equation}
On the other hand,
from (\ref{eq3}),
there exists a bi-homogeneous polynomial $U_0$ of bi-degree $(a-q-1, b)$
such that
\begin{equation}\label{eq5}
U(X_0, X_1;Y_0, Y_1) = U_0(X_0, X_1;Y_0, Y_1)X_0 + f_1(Y_0, Y_1)X_1^{a-q}.
\end{equation}
Let 
\[
V_0 = X_1^{a-q} f_2(Y_0, Y_1)(X_0^{q-1}X_1 - X_1^q) + V,
\]
which is of bi-degree $(a, b-q)$.
Then, from (\ref{eq1}),  (\ref{eq4}) and  (\ref{eq5}),
\begin{equation}\label{eq6}
   \begin{split}
F(X_0, X_1;Y_0, Y_1)=&U_0(X_0, X_1;Y_0, Y_1) (X_0^qX_1 - X_0X_1^q) \\
& + V_0(X_0, X_1;Y_0, Y_1) (Y_0^{q-1}Y_1 -Y_1^q).
   \end{split}
\end{equation}
Substitute $Y_0=0$ in (\ref{eq6}).
Then we have
\[
F(X_0, X_1;0, Y_1)=U_0(X_0, X_1;0, Y_1) (X_0^qX_1 - X_0X_1^q) 
- V_0(X_0, X_1;0, Y_1) Y_1^q.
\]
Since $V_0$ is bi-homogeneous of degree $(a, b-q)$,
there is a homogeneous polynomial $g_1(X_0,X_1)$ of degree $a$
such that $V_0(X_0, X_1;0, Y_1)=g_1(X_0, X_1) Y_1^{b-q}.$
Since $F(\alpha_0, \alpha_1;0, 1)= 0$ for any
$(\alpha_0, \alpha_1) \in \mathbb{P}^1(\mathbb{F}_q)$
and the first term in the right side of (\ref{eq6})
vanishes at $(\alpha_0, \alpha_1; 0,1)$
for any $(\alpha_0, \alpha_1) \in \mathbb{P}^1(\mathbb{F}_q)$,
so is $g_1(X_0, X_1)$.
Hence
\[
g_1(X_0, X_1) = f_3(X_0, X_1) (X_0^qX_1 - X_0X_1^q)
\]
for a homogeneous polynomial $ f_3(X_0, X_1)$
of degree $a-(q+1)$.
Therefore there is a bi-homogeneous polynomial $g(X_0, X_1;Y_0, Y_1)$
of bi-degree $(a, b-q-1)$
such that
\begin{equation}\label{eq7}
  \begin{split}
 V_0(X_0, X_1;Y_0, Y_1) =& g(X_0, X_1;Y_0, Y_1)Y_0\\
   & +f_3(X_0,X_1)Y_1^{b-(q+1)}(X_0^qX_1 - X_0X_1^q).
  \end{split}
\end{equation}
Finally, taking $f(X_0, X_1;Y_0, Y_1)$ as
\[
f = U_0 +f_3(X_0, X_1) Y_1^{b-(q+1)},
\]
we have
\[
F = f \cdot (X_0^q X_1 - X_0 X_1^q) + g \cdot(Y_0^q Y_1 - Y_0 Y_1^q)
\]
by (\ref{eq6}) and (\ref{eq7}).
\qed

\section{Existence for $q\neq 2$}
Our method for proving the existence is rather heuristic.
In general, for a curve $C$ on $\mathbb{P}^1\times \mathbb{P}^1$
of bi-degree $(a,b)$ defined by
$F=F(X_0, X_1;Y_0, Y_1)=0$,
if $C$ is nonsingular with $a>0$ and $ b>0$, then $C$ is irreducible
\footnote{Actually, if $C$ is reducible, we can pick up two components $D_1$ of type $(a', b')$ with $a' >0$ and $D_2$ of type $(a'', b'')$ with $b''>0$.
Then $(D_1.D_2) = a'b''+a''b' >0$,
which implies $D_1\cap D_2 \neq \emptyset$.};
and the nonsingularity of $C$ is guaranteed by the condition that
the simultaneous equations
\begin{equation}\label{eq3-1}
F=F_{X_0}=F_{X_1}=F_{Y_0}=F_{Y_1}=0
\end{equation}
has only trivial solutions.
Here the subscript variable of a polynomial indicates
taking the partial derivative of the polynomial by the assigned variable,
and a trivial solution means a solution of the form
$(0,0)\times (\ast, \ast)$ or $(\ast, \ast)\times (0,0) $.
The strategy to show the existence is to find a polynomial $F$ equipped with the above property.

In this section, we suppose that $q \neq 2$,
and restrict our consideration within the following situation.
\begin{setup}\label{setup}
Let $C$ be a curve on $\mathbb{P}^1\times \mathbb{P}^1$ over $\mathbb{F}_q$
of bi-degree $(q+1, q+1)$
defined by $F=0$.
Here
$F=F(X_0, X_1;Y_0, Y_1)$
is defined as
\[
F= f(Y_0, Y_1) 
\left(X_0^qX_1 -X_0X_1^q\right) +
g(X_0, X_1) 
\left(Y_0^qY_1 -Y_0Y_1^q\right),
\]
where $f(Y_0, Y_1)$ and  $g(X_0, X_1)$
are homogeneous polynomials of degree $q+1$.
Furthermore, suppose that
both $f$ and $g$ have no multiple zeros on $\mathbb{P}^1$,
and no $\mathbb{F}_q$-zeros.
\end{setup}

\begin{lemma}\label{keylemma}
Under Setup~{\rm \ref{setup}},
a point on $\mathbb{P}^1 \times \mathbb{P}^1$ satisfies 
the equations {\rm (\ref{eq3-1})} if and only if it satisfies
\begin{equation}\label{eq3-2}
\frac{g_{X_1}}{X_0^q} = - \frac{g_{X_0}}{X_1^q}
=-\frac{f_{Y_1}}{Y_0^q} = \frac{f_{Y_0}}{Y_1^q}.
\end{equation}
\end{lemma}
\proof
Let us consider five polynomials
\begin{align}
F =& f(Y_0,Y_1)(X_0^qX_1- X_0X_1^q) + g(X_0, X_1)(Y_0^qY_1- Y_0Y_1^q) 
                                       \label{eqF1} \\
F_{X_0} =& -X_1^q  f(Y_0,Y_1) + g_{X_0}(X_0, X_1) (Y_0^qY_1- Y_0Y_1^q)
                                       \label{eqF2} \\
F_{X_1} =& X_0^q  f(Y_0,Y_1) + g_{X_1}(X_0, X_1)(Y_0^qY_1- Y_0Y_1^q)
                                       \label{eqF3} \\
F_{Y_0} =& f_{Y_0}(Y_0,Y_1) (X_0^qX_1- X_0X_1^q) - Y_1^q g(X_0, X_1)
                                       \label{eqF4} \\
F_{Y_1} =& f_{Y_1}(Y_0,Y_1) (X_0^qX_1- X_0X_1^q) + Y_0^q g(X_0, X_1).
                                       \label{eqF5}
\end{align}

\noindent
{\sl Step}~1.
Let
$(\alpha_0, \alpha_1) \times (\beta_0, \beta_1) 
\in \mathbb{P}^1 \times \mathbb{P}^1$
be a solution of (\ref{eq3-1}).
Then
\begin{enumerate}[(i)]
\item $(\alpha_0, \alpha_1) \not\in \mathbb{P}^1(\mathbb{F}_q)$
and  $(\beta_0, \beta_1) \not\in \mathbb{P}^1(\mathbb{F}_q)$;
\item $f(\beta_0, \beta_1) \neq 0$ and $g(\alpha_0, \alpha_1)\neq 0$;
\item $g_{X_0}(\alpha_0, \alpha_1)\neq 0$,
      $g_{X_1}(\alpha_0, \alpha_1)\neq 0$,
$f_{Y_0}(\beta_0, \beta_1)\neq 0$ 
         and $f_{Y_1}(\beta_0, \beta_1)\neq 0$.
\end{enumerate}

For (i), suppose that $(\alpha_0, \alpha_1) \in \mathbb{P}^1(\mathbb{F}_q)$.
Since $(\alpha_0^q\alpha_1 - \alpha_0\alpha_1^q)=0$ and 
$F_{Y_0}(\alpha_0, \alpha_1; \beta_0, \beta_1)= 
F_{Y_1}(\alpha_0, \alpha_1; \beta_0, \beta_1)=0$,
we have $g(\alpha_0, \alpha_1)=0$ by
(\ref{eqF4}) and (\ref{eqF5}), which contradicts the setup on $g$.
Similarly, we have  $(\beta_0, \beta_1) \not\in \mathbb{P}^1(\mathbb{F}_q)$.

For (ii), 
suppose that $f(\beta_0, \beta_1)=0$.
Then
\[
0= F(\alpha_0, \alpha_1; \beta_0, \beta_1) 
= g(\alpha_0, \alpha_1)(\beta_0^q \beta_1 - \beta_0\beta_1^q).
\]
Since $(\beta_0, \beta_1) \not\in \mathbb{P}^1(\mathbb{F}_q)$ by (i),
$g(\alpha_0, \alpha_1)=0.$
Similarly, we have $g_{X_0}(\alpha_0, \alpha_1)=0$ and
$g_{X_1}(\alpha_0, \alpha_1)=0$ by using (\ref{eqF2}) and (\ref{eqF3}).
So $(\alpha_0, \alpha_1)$ is a multiple zero of $g$, which contradicts
the setup on $g$.

For (iii),
suppose that, for example, $g_{X_0}(\alpha_0, \alpha_1)=0$.
Then, from (\ref{eqF2}) and a result in (ii),
$\alpha_1=0$.
Namely $(\alpha_0, \alpha_1)=(1,0)$, which contradicts
a result in (i).

\noindent
{\sl Step}~2.
Continuously, 
let $(\alpha_0, \alpha_1) \times (\beta_0, \beta_1) $
be a solution of (\ref{eq3-1}).
Then
\begin{equation}\label{eqF6}
\beta_0^q \beta_1 - \beta_0\beta_1^q = 
\frac{\alpha_1^q f(\beta_0, \beta_1)}{g_{X_0}(\alpha_0, \alpha_1)}
=
- \frac{\alpha_0^q f(\beta_0, \beta_1)}{g_{X_1}(\alpha_0, \alpha_1)}
\end{equation}
by (\ref{eqF2}) and (\ref{eqF3}).
Similarly,
\begin{equation}\label{eqF7}
\alpha_0^q \alpha_1 - \alpha_0\alpha_1^q = 
\frac{\beta_1^q g(\alpha_0, \alpha_1)}{f_{Y_0}(\beta_0, \beta_1)}
=
- \frac{\beta_0^q g(\alpha_0, \alpha_1)}{f_{Y_1}(\beta_0, \beta_1)}
\end{equation}
by (\ref{eqF6}) and (\ref{eqF5}).
Since $ f(\beta_0, \beta_1) \neq 0$ and
$ g(\alpha_0, \alpha_1) \neq 0$
by Step~1~(ii),
(\ref{eqF6}) and (\ref{eqF7})
imply that the point $(\alpha_0, \alpha_1) \times (\beta_0, \beta_1) $
satisfies the equations
\[
\frac{X_1^q}{g_{X_0}} = - \frac{X_0^q}{g_{X_1}}
\mbox{\quad  and \quad }
 \frac{Y_1^q}{f_{Y_0}}= -\frac{Y_0^q}{f_{Y_1}}.
\]
On the other hand, substitute the first equality of (\ref{eqF6}) and 
that of (\ref{eqF7}) into
\[
f(\beta_0, \beta_1)(\alpha_0^q \alpha_1 - \alpha_0\alpha_1^q)
+ g(\alpha_0, \alpha_1)(\beta_0^q \beta_1 - \beta_0\beta_1^q)=0,
\]
we know that the point also satisfies
$\frac{X_1^q}{g_{X_0}} = - \frac{Y_1^q}{f_{Y_0}}.$
So the point satisfies (\ref{eq3-2}).

\noindent
{\sl Step}~3.
Conversely, suppose that (\ref{eq3-2}) holds at a point
$(\alpha_0, \alpha_1) \times (\beta_0, \beta_1). $
Those equations imply the following four equations:
\begin{align}
X_0^q f_{Y_1} + Y_0^q g_{X_1} &=0  \label{4eq1}\\
X_0^q f_{Y_0} - Y_1^q g_{X_1} &=0  \label{4eq2}\\
X_1^q f_{Y_1} - Y_0^q g_{X_0} &=0  \label{4eq3}\\
X_1^q f_{Y_0} + Y_1^q g_{X_0} &=0.  \label{4eq4}
\end{align}
Making
\[
[\mbox{Eq.~(\ref{4eq1})}] \times X_1Y_1
+[\mbox{Eq.~(\ref{4eq2})}] \times X_1Y_0
-[\mbox{Eq.~(\ref{4eq3})}] \times X_0Y_1
-[\mbox{Eq.~(\ref{4eq4})}] \times X_0Y_0,
\]
we have
\[
(Y_0f_{Y_0} +Y_1f_{Y_1})(X_0^qX_1-X_0X_1^q)
+ (X_0g_{X_0} +X_1g_{X_1})(Y_0^qY_1-Y_0Y_1^q)=0,
\]
which means $F=0$ by the Euler identity.
Similarly,
making
$[\mbox{Eq.~(\ref{4eq1})}] \times Y_1
+[\mbox{Eq.~(\ref{4eq2})}] \times Y_0,$
$-\left([\mbox{Eq.~(\ref{4eq3})}] \times Y_1
+[\mbox{Eq.~(\ref{4eq4})}] \times Y_0\right),$
$[\mbox{Eq.~(\ref{4eq1})}] \times X_1
-[\mbox{Eq.~(\ref{4eq3})}] \times X_0$
and
$[\mbox{Eq.~(\ref{4eq2})}] \times X_1
-[\mbox{Eq.~(\ref{4eq4})}] \times X_0$,
we have
$F_{X_1}=0$,
$F_{X_0}=0$,
$F_{Y_1}=0$
and
$F_{Y_0}=0$
respectively.
\qed

\begin{theorem}
If $q \neq 2$, then there is a nonsingular irreducible
$\mathbb{P}^1 \times \mathbb{P}^1$-filling curve
over $\mathbb{F}_q$ of bi-degree $(q+1, q+1)$.
\end{theorem}
\proof
From Lemma~\ref{keylemma},
it is sufficient to find two homogeneous polynomials
$f(Y_0, Y_1)$ and $g(X_0, X_1)$ such that
\begin{enumerate}[(i)]
\item $f=0$ and $g=0$ have no multiple solutions on $\mathbb{P}^1$;
\item $f=0$ and $g=0$ have no $\mathbb{F}_q$-solutions:
\item the simultaneous equations (\ref{eq3-2})
has no nontrivial solutions.
\end{enumerate}

\noindent
(I) Suppose that $q$ is odd and $q\geq 5$.
Since
\[
|\mathbb{F}_q^{\ast} \setminus \{u^2 \mid u \in \mathbb{F}_q^{\ast}\}|
 = \frac{q-1}{2} \geq 2, 
\]
we can choose two distinct elements
$-\delta, -\gamma \in 
\mathbb{F}_q^{\ast} \setminus \{u^2 \mid u \in \mathbb{F}_q^{\ast}\}$.
Let
$f(Y_0, Y_1) = Y_0^{q+1} + \delta Y_1^{q+1}$
and
$g(X_0, X_1) = X_0^{q+1} + \gamma X_1^{q+1}.$
Obviously, they satisfy the condition (i).
By the choice of $\delta$ and $\gamma$, they satisfy the condition (ii).
For our $f$ and $g$ above,
the simultaneous equations (\ref{eq3-2}) can be reduced to
\[
\frac{\gamma X_1^q}{X_0^q} = -\frac{X_0^q}{X_1^q}= - \frac{\delta Y_1^q}{Y_0^q}
=\frac{Y_0^q}{Y_1^q}.
\]
Put $\frac{X_1}{X_0} =x$ and $\frac{Y_1}{Y_0} =y$.
Then
\[
\gamma x^q = -\frac{1}{x^q} = -\delta y^q = \frac{1}{y^q},
\]
which imply
\[
\gamma =  -\frac{1}{x^{2q}},\   \delta = -\frac{1}{y^{2q}} 
\mbox{\quad and \quad} x^q = -y^q.
\]
Hence $\gamma = \delta$,
which is a contradiction.

\noindent
(II) Suppose that $q=2^e$ with $e \geq 2$.
Since
\[
|\mathbb{F}_{2^e} \setminus \{u + u^2 \mid u \in \mathbb{F}_{2^e}\}|
 = 2^{e-1} \geq 2,
\]
we can choose two distinct elements
$\delta , \gamma \in
\mathbb{F}_{2^e} \setminus \{u + u^2 \mid u \in \mathbb{F}_{2^e}\}.$
Let
$f(Y_0, Y_1) = Y_0^{q+1} + Y_0 Y_1^q + \delta Y_1^{q+1}$
and
$g(X_0, X_1) = X_0^{q+1} + X_0 X_1^q + \gamma X_1^{q+1}.$
Then, by similar arguments in (I), we know that
$f$ and $g$ satisfy the above conditions (i), (ii) and (iii).

\noindent
(III) Suppose that $q=3$.
Let
$
f(Y_0, Y_1) = Y_0^4 + Y_1^4
$
and
$
g(X_0, X_1) = X_0^4 + X_0X_1^3 +2X_1^4.
$
It is easy to see that $f$ and $g$ satisfy the conditions (i) and (ii).
For our $f$ and $g$,
the equations (\ref{eq3-2}) can be reduced to
\[
\frac{2 X_1^3}{X_0^3} =2 \frac{X_0^3+ X_1^3}{X_1^3}= 2\frac{Y_1^3}{Y_0^3}= \frac{Y_0^3}{Y_1^3}.
\]
Put $\frac{X_1}{X_0} =x$ and $\frac{Y_1}{Y_0} =y$.
Then
the simultaneous equations can be written as
\[
2x^3 = 2\left( \frac{1}{x^3} +1 \right) = 2 y^3 = \frac{1}{y^3},
\]
and has no solutions.
Actually,
from the first and third terms, we have $x=y$.
Then the second and the forth terms imply
$x=y=1$.
But this solutions do not agree with the first equality.
\qed

\section{Non-existence and existence for $q = 2$}
\begin{theorem}
Assume that $q=2$.
\begin{enumerate}[{\rm (a)}]
\item There are no irreducible $\mathbb{P}^1 \times \mathbb{P}^1$-filling curves over $\mathbb{F}_2$ of bi-degree $(3,3)$.
\item There are nonsingular irreducible $\mathbb{P}^1 \times \mathbb{P}^1$-filling curves over $\mathbb{F}_2$ of bi-degree $(4,3)$, and of bi-degree $(3,4)$.
\end{enumerate}
\end{theorem}
\proof
(a) Let $C$ be an irreducible curve on $\mathbb{P}^1 \times \mathbb{P}^1$
over $\mathbb{F}_2$ of bi-degree $(3,3)$.
Through the Segre embedding
$s : \mathbb{P}^1 \times \mathbb{P}^1 \hookrightarrow \mathbb{P}^3$,
$s(C)$ is a nondegenerate irreducible curve of degree $6$.
Therefore, by \cite[Theorem~3.2]{hom2012},
\[
N_2(C) \leq
\left\lfloor\frac{(2-1)(2^{4}-1)}{2(2^3-1) -3(2-1)}6 \right\rfloor
= \left\lfloor\frac{90}{11}\right\rfloor < 9
=N_2(\mathbb{P}^1 \times \mathbb{P}^1).
\]

(b) Let us consider the curve $C$ on $\mathbb{P}^1 \times \mathbb{P}^1$
over $\mathbb{F}_2$ of bi-degree $(4,3)$ defined by
$F=0$, where
\[
F=
(X_0Y_0^3 + X_1Y_1^3) (X_0^2 X_1 + X_0 X_1^2)
   + (X_0^2 + X_0 X_1 +X_1^2)^2 (Y_0^2 Y_1 + Y_0 Y_1^2).
\]

We want to confirm that $C$ is nonsingular.
Partial derivatives of $F$ by each variable are
\begin{align*}
F_{X_0} =& Y_0^3(X_0^2 X_1 + X_0 X_1^2)+
                  X_1^2  (X_0Y_0^3 + X_1Y_1^3) \\
F_{X_1} =&   Y_1^3(X_0^2 X_1 + X_0 X_1^2)+
                  X_0^2  (X_0Y_0^3 + X_1Y_1^3) \\
F_{Y_0} =& X_0 Y_0^2 (X_0^2 X_1 + X_0 X_1^2)+
                 (X_0^2 + X_0 X_1 +X_1^2)^2Y_1^2 \\
F_{Y_1} =&   X_1 Y_1^2 (X_0^2 X_1 + X_0 X_1^2)+
                 (X_0^2 + X_0 X_1 +X_1^2)^2Y_0^2 .
\end{align*}
We will try to find a nontrivial solution for (\ref{eq3-1}),
and expect this trial to fail.

If $X_0=0$, then $X_1^4Y_1^2=0$ and $X_1^4Y_0^2=0$
from $F_{Y_0} =0$ and $F_{Y_1} =0$ respectively, which means
that we get only trivial solutions.

If $Y_0 =0$, then $X_1^3 Y_1^3=0$ by $F_{X_0}=0$.
Since we are seeking a nontrivial solution,
$X_1 =0$. Hence $X_0^4Y_1^2=0$ by $F_{Y_0}=0$,
which is absurd.

Let $x = \frac{X_1}{X_0}$
and $y = \frac{Y_1}{Y_0}$.
Then the problem can be reduced to solve
the simultaneous equations:
\begin{align}
 & (1+xy^3)(x+x^2) + (1+x+x^2)^2(y+y^2) =0  \label{nh2F} \\
 & (x+x^2) + x^2(1+xy^3) =0 \label{nh2FX0} \\
 &  (x+x^2)y^3 + (1+xy^3) =0 \label{nh2FX1} \\
 & (x+x^2) + (1+x+x^2)^2 y^2=0 \label{nh2FY0}\\
 & x(x+x^2)y^2 + (1+x+x^2)^2 =0. \label{nh2FY1}
\end{align}

Suppose that $(x,y)$ is a solution of the above simultaneous equations.
Then $x \neq 0, 1$. Actually $x=0$ is not compatible with (\ref{nh2FX1}),
neither is $x=1$ with (\ref{nh2FY1}).

Making $\text{(\ref{nh2FX0})}+ x^2\times \text{(\ref{nh2FX1})}$
we get
\[
(x+x^2)(1+x^2y^3) =0.
\]
Hence $x^2y^3=1$. On the other hand,
making $\text{(\ref{nh2FY0})}+ y^2\times \text{(\ref{nh2FY1})}$,
\[
 (x+x^2)(1 + x y^4)=0.
\]
Hence $xy^4 =1$.
Therefore $x =y \neq 0$.
Substitute $y=x$ into (\ref{nh2FX0}) and (\ref{nh2FY0}).
Then we get $x+x^6=0$ and $x + x^4 +x^6 =0$,
which is absurd.
\qed

\section{Closing}
The image of
the Segre embedding 
$
\mathbb{P}^1 \times \mathbb{P}^1 \hookrightarrow \mathbb{P}^3
$
over $\mathbb{F}_q$
is a nonsingular quadratic surface,
and is called  a hyperbolic quadratic surface \cite[5.2]{hir}.
\begin{corollary}\label{finalcor}
Let $S$ be a hyperbolic quadratic surface over $\mathbb{F}_q$
in $\mathbb{P}^3$.
The smallest degree of an $S$-filling irreducible curve is
$2q+2$ if $q\neq 2$, and $7$ if $q=2$.
\end{corollary}

\begin{remark}
We have constructed a curve in Corollary~\ref{finalcor}
as a nonsingular one.
But there is a nonsingular $S$-filling curve which is not irreducible;
the curve $\cup_{P \in \mathbb{P}^1(\mathbb{F}_q)} \mathbb{P}^1\times \{P \}$
on $\mathbb{P}^1 \times \mathbb{P}^1 \simeq S$
is such a curve, and of degree $q+1$.
\end{remark}

At the end of this paper,
we pose a problem.
For a given nonsingular surface $S$ over $\mathbb{F}_q$ in $\mathbb{P}^3$,
the existence of an $S$-filling nonsingular curve is known (\cite{gab} and \cite{poo}), however, we don't know whether an $S$-filling irreducible curve exists or not.
But we strongly expect the existence of an $S$-filling irreducible curve.
A problem is what the smallest degree of an $S$-filling irreducible curve is,
if any.
Especially the following three cases may be interesting:
\begin{enumerate}[(i)]
\item the elliptic quadratic surface,
\item the Hermitian surface of degree $\sqrt{q}+1$,
\item the space filling surface 
\[
X_0X_1^q - X_0^qX_1 + X_2X_3^q - X_2^qX_3=0.
\]
\end{enumerate}

\end{document}